\documentclass[11pt]{article}
\usepackage{amsmath,amsfonts, amssymb,amsthm}

\bibstyle{plain}

\def\be{\begin{equation}}
\def\en{\end{equation}}
\def\bee{\begin{eqnarray*}}
\def\ene{\end{eqnarray*}}

\def\E{{\bf E}}
\def\P{{\mathbb{P}}}
\def\R{{\mathbb{R}}}
\def\E{{\mathbb{E}}}

\def\ep{\varepsilon}

\vspace{1in}

\title{\protect \large \bf
CENTRAL LIMIT THEOREM AND \\ 
DIOPHANTINE APPROXIMATIONS
\thanks{Partially supported by the NSF grant DMS-1612961}
}

\vspace{2in}

\author{Sergey G. Bobkov
\footnote{Address: School of Mathematics, University of Minnesota, 
127 Vincent Hall, 206 Church St. S.E., \newline
\hskip10mm Minneapolis, MN 55455 USA. \ Email: bobkov@math.umn.edu}
}

\begin{document}

\bibliographystyle{plain}

\maketitle

\begin{abstract}
\hskip-6mm
Let $F_n$ denote the distribution function of the normalized sum
$Z_n = (X_1 + \dots + X_n)/\sigma\sqrt{n}$ of i.i.d. random variables
with finite fourth absolute moment. In this paper, polynomial rates of 
convergence of $F_n$ to the normal law with respect to the Kolmogorov 
distance, as well as polynomial approximations of $F_n$ by the Edgeworth 
corrections (modulo logarithmically growing factors in $n$) are
given in terms of the characteristic function of $X_1$. Particular cases 
of the problem are discussed in connection with Diophantine approximations.
\end{abstract}

\section{Introduction}
\setcounter{equation}{0}

Let $X,X_1,X_2,\dots$ be independent, identically distributed
random variables with mean zero, variance $\sigma^2$ $(\sigma > 0)$
and finite 3-rd absolute moment $\beta_3 = \E\,|X|^3$. 
Denote by $F(x) = \P\{X \leq x\}$ the distribution function
and by $f(t) = \E\,e^{itX}$ the characteristic function of $X$.

The Berry-Esseen theorem provides a standard rate of approximation of
the distribution functions $F_n(x) = \P\{Z_n \leq x\}$ of the normalized sums
$$
Z_n = \frac{X_1 + \dots + X_n}{\sigma \sqrt{n}}
$$
by the standard normal distribution function $\Phi(x)$ with density
$\varphi(x) = \frac{1}{\sqrt{2\pi}}\,e^{-x^2/2}$ ($x \in \R$).
Namely, up to a numerical constant $c$, we have
$$
\sup_x |F_n(x) - \Phi(x)| \leq c\,
\frac{\beta_3}{\sigma^3\sqrt{n}}.
$$
In general, higher order moment assumptions do not improve this rate, 
as can been seen on the example of lattice distributions $F$.
Nevertheless, under the Cram\'er condition
\be
\limsup_{t \rightarrow \infty}\, |f(t)| < 1,
\en
it is possible to slightly correct the limit law (by allowing dependence
in $n$), so as to improve the rate of approximation. In particular, 
consider an Edgeworth correction of the 3-rd order
\be
\Phi_3(x) = \Phi(x) - 
\frac{\alpha_3}{6\sigma^3\sqrt{n}}\,(x^2-1)\, \varphi(x), \qquad 
\alpha_3 = \E X^3,
\en
which also depends on $n$, except for the case $\alpha_3 = 0$ 
(when $\Phi_3 = \Phi$). It is well-known that, if the 4-th absolute moment 
$\beta_4 = \E X^4$ is finite, the uniform deviations 
$$
\Delta_n = \sup_x |F_n(x) - \Phi_3(x)|
$$
are at most of order $1/n$.
Moreover, with higher order moments assumptions, the corresponding higher
order Edgeworth corrections (called also Edgeworth expansions) 
provide an error of approximation
decaying as powers of $1/\sqrt{n}$, cf. e.g. [P1], [B-RR].

Without the Cram\'er condition (1.1), the problem of possible rates is 
rather delicate, as the order of magnitude of $\Delta_n$ depends on 
arithmetical properties of the point spectrum of $F_n$. This was
already emphasized by Esseen, who established the following general result
(cf. [E], pp. 49-53): If $X$ has a non-lattice distribution (equivalently, 
$|f(t)|<1$ for all $t > 0$), and if the 3-rd absolute moment 
of $X$ is finite, then
\be
\Delta_n = o\Big(\frac{1}{\sqrt{n}}\Big) \quad {\rm as} \ \ 
n \rightarrow \infty.
\en
It seems that not much has been said in literature in addition to this theorem 
(see, however, a cycle of papers [Ch]).
The aim of these notes is to refine (1.3) by connecting
possible polynomial rates for $\Delta_n$ with behavior of the characteristic 
function $f(t)$ at infinity. Let us stress that, although
the lack of the Cram\'er property forces $F$ not to have an absolutely
continuous component, the class of probability distributions with
$\limsup_{t \rightarrow \infty}\, |f(t)| = 1$ is extremely rich and
interesting
(including discrete and many purely singular continuous probability measures).

For simplicity, we focus on intermediate rates between 
$\frac{1}{\sqrt{n}}$ and $\frac{1}{n}$ for $\Delta_n$. Let us state the
relationship, by using the notation
$\widetilde O(t^p)$ for the growth rate $O(t^p\,(\log t)^q)$ with some $q \in \R$, 
and similarly $\widetilde O(n^{-p})$ for $O(n^{-p}\,(\log n)^q)$.

\vskip5mm
{\bf Theorem 1.1.} {\sl Suppose that $\beta_4 < \infty$. Given 
$p \geq 2$, the following two properties are equivalent:
\begin{eqnarray}
\frac{1}{1 - |f(t)|} 
 & = &
\widetilde O(t^p\,) \qquad \quad {\sl as} \ \  t \rightarrow \infty; \\
\Delta_n
 & = &
\widetilde O\big(n^{-\frac{1}{2} - \frac{1}{p}}\big) \quad {\sl as} \ \ 
n \rightarrow \infty.
\end{eqnarray}
}

A more precise formulation reflecting appearance of the logarithmic
factors in $\widetilde O$ in (1.4)-(1.5) will be given in Sections 3 
and 5. As for the restriction $p \geq 2$, it may actually be relaxed to $p>0$
under higher moment assumptions by adding to $\Phi_3$ other terms in 
the corresponding Edgeworth expansions. 

Let us illustrate Theorem 1.1 in a simple discrete situation.
As is standard, we denote by 
$\|x\|$ the distance from a real number $x$ to the closest integer.
Given an irrational real number $\alpha$, define the quantity
$$
\eta(\alpha) \, = \, \sup\Big\{\eta > 0:
\liminf_{n \rightarrow \infty}\, n^{\eta} \|n\alpha\| = 0\Big\} \, = \,
\inf\Big\{\eta > 0:\, \inf_{n \geq 1}\, n^{\eta} \|n\alpha\| > 0\Big\}.
$$
One says that $\alpha$ is of type $\eta = \eta(\alpha)$ and calls
$1 + \eta$ an irrationality exponent of $\alpha$. Equivalently, the value 
of $\eta$ is an optimal one, for which the Diophantine inequality
$$
\Big|\alpha - \frac{p}{q}\Big| < \frac{1}{q^{1 + \eta - \ep}}
$$
has infinitely many rational solutions $\frac{p}{q}$ with any fixed $\ep>0$ 
(cf. e.g. [K-N], [B-B-S]). Thus, this quantity provides an important
information on how well the number $\alpha$ may be approximated by rationals.
By Dirichlet's theorem, necessarily $\eta \geq 1$, and actually the possible values 
of $\eta$ fill the half-axis $[1,\infty]$ including the case $\eta = \infty$ 
(which describes the class of Liouville's numbers).

Applying Theorem 1.1 with $p = 2\eta$, one may derive the next 
characterization.

\vskip5mm
{\bf Corollary 1.2.} {\sl Given an irrational number $\alpha$, suppose that 
the random variable $X$ takes the values $\pm 1$ and $\pm \alpha$ each with 
probability $1/4$. Then $\alpha$ is of finite type $\eta$, if and only if, for any $\ep > 0$,
\be
\sup_x\, |F_n(x) - \Phi(x)| \, = \,
O\left(n^{-\frac{1}{2} - \frac{1}{2\eta} + \ep}\right) \quad {\sl as} \ \ 
n \rightarrow \infty.
\en
}

\vskip2mm
A similar description continuous to hold when $X$ takes the values 
$\pm 1 \pm \alpha$. In this case, one may write $X = X' + \alpha X''$ 
in the sense of laws, 
where $X'$ and $X''$ are independent random variables with a symmetric 
Bernoulli distribution on $\{-1,1\}$. While for $X'$ and $\alpha X''$
separately, the corresponding deviations $\Delta_n$ are of order $1/\sqrt{n}$,
we see that the convolution structure in the underlying distribution $F$ 
may essentially improve the rate.

For example, by Roth's theorem (cf. [C], [S1-2]), we have 
$\eta = 1$ for any irrational algebraic $\alpha$, and then (1.6) becomes
$\Delta_n = O(n^{-1 + \ep})$. If $\alpha$ is a quadratic irrationality,
or more generally, a badly approximable number, one may sharpen the rate
to $\Delta_n = O(\frac{1}{n}\sqrt{\log n}\,)$.
Although in these examples, such $\alpha$'s form a set
of (Lebesgue) measure zero, a slightly worse rate
$$
\Delta_n = O\Big(\frac{1}{n}\,(\log n)^{\frac{3}{2} + \ep}\Big)
$$
can be derived for almost all values of $\alpha$ on the line (see Section 7 for details).

It is interesting to compare relation (1.6) with a statement about
an asymptotic behavior of ``empirical" measures
$$
\widetilde F_n = \frac{1}{n}\,\sum_{k=1}^n \delta_{\{k\alpha\}},
$$
where $\{x\}$ stands for the fractional part and $\delta_x$ denotes
a point mass at a given point (one may similarly consider the  sequence
$\|k\alpha\|$ and use the identity $\|x\| = \min\{\{x\},1-\{x\}\}$).
By Weyl's criterion,
$\widetilde F_n$ are weakly convergent to the uniform distribution
on $(0,1)$, as long as $\alpha$ is irrational. Results by
Hecke, Ostrowski and Behnke in 1920's quantify this convergence:
For any $\ep>0$, with some positive $c_0 = c_0(\alpha,\ep)$ and 
$c_1 = c_1(\alpha,\ep)$, we have
\be
c_0\,n^{-\frac{1}{\eta} - \ep} \, \leq \,
\sup_{0<x<1} \big|\widetilde F_n(x) - x\big| \, \leq \, 
c_1\,n^{-\frac{1}{\eta} + \ep},
\en
where $\eta = \eta(\alpha)$ ([K-N]).
Although there is some difference between (1.6) and (1.7), the two rates
turn out to be in essence the same in the critical case $\eta=1$.
Let us also mention that, for quadratic irrationalities $\alpha$,
an asymptotic behavior of $\widetilde F_n$ has been comprehensively
studied in the recent times by Beck [Be].

The paper is organized as follows. In section 2 we remind a basic
Berry-Esseen-type bound for the distributions $F_n$
which is applicable to reach the rate of approximation of $F_n$ by $\Phi_3$
potentially up to order $1/n$. Here we also explain the sufficiency part 
in Theorem 1.1. In Sections 3-4 we discuss non-uniform
bounds on $|F_n(x) - \Phi_3(x)|$ together with bounds on the difference
between the Fourier-Stieltjes transforms of $F_n$ and $\Phi_3$.
The necessity part in Theorem 1.1 is considered separately in Section 5.
Section 6 deals with Diophantine inequalities, where Corollary 1.2 is 
derived, actually in a somewhat more general and precise form.
Applications of this corollary are clarified in Section 7.

\vskip10mm
\section{Berry-Esseen inequality. Sufficiency part in Theorem 1.1}
\setcounter{equation}{0}

The derivation of uniform estimates on the difference between distribution 
functions, say $F$ and $G$, is commonly based on a general Berry-Esseen bound 
\be
c\,\sup_x\,|F(x) - G(x)| \, \leq \,
\int_0^T \frac{|f(t) - g(t)|}{t}\,dt + \frac{D}{T} \qquad (T>0),
\en
involving the Fourier-Stieltjes transforms
$$
f(t) = \int_{-\infty}^\infty e^{itx}\,dF(x), \quad
g(t) = \int_{-\infty}^\infty e^{itx}\,dG(x) \qquad (t \in \R).
$$
Here and below we denote by $c$ a positive absolute 
constant which may be different in different places.
In fact, in (2.1), $G$ may be an arbitrary differentiable function 
of bounded variation on the real line such that $G(-\infty)=0$, $G(\infty)=1$,
and $\sup_x\,|G'(x)| \leq D$ (cf. [E], [P2], [Bo1]).
With this approach, the implication $(1.4) \Rightarrow (1.5)$ is rather
standard (although we cannot give an exact reference). For completeness, 
we remind the basic argument in the
special situation as in Theorem~1.1 which yields an upper bound on the
uniform distance
$$
\Delta_n = \sup_x\, |F_n(x) - \Phi_3(x)|.
$$

Namely, one may apply (2.1) with
$F_n$ in place of $F$ and with $G = \Phi_3$. The Fourier-Stieltjes 
transform of $F_n$ is just the characteristic function of $Z_n$ given by 
$f_n(t) = f(\frac{t}{\sigma \sqrt{n}})^n$, where $f$ is 
the characteristic function of $X$. The Fourier-Stieltjes 
transform of $\Phi_3$ is 
\be
g_3(t) = e^{-t^2/2} +
\frac{\alpha_3}{6\sigma^3\sqrt{n}}\, (it)^3\,e^{-t^2/2} \qquad (t \in \R).
\en
Such an application then leads to the following estimate.

\vskip5mm
{\bf Lemma 2.1.} {\sl Suppose that $\beta_4$ is finite. 
For all $n \geq 1$ and $T \geq \frac{\sigma}{\sqrt{\beta_4}}$,
\be
c\,\Delta_n \, \leq \, \frac{\beta_4}{\sigma^4 n} + 
\frac{1}{T\sigma \sqrt{n}} + 
\int_{\frac{\sigma}{\sqrt{\beta_4}}}^T \frac{|f(t)|^n}{t}\,dt.
\en
}

\vskip2mm
{\bf Proof.} Put $T_0 = \frac{\sigma^2}{\sqrt{\beta_4}} \sqrt{n}$ and
introduce the Lyapunov coefficients 
$L_s = \frac{\beta_s}{\sigma^s}\,n^{-\frac{s-2}{2}}$ ($\beta_s = \E\,|X|^s$),
which we need for $s=3$ and $s=4$. Since the function 
$s \rightarrow L_s^{1/(s-2)}$ is non-decreasing in $s > 2$, we have 
$L_3 \leq L_4^{1/2}$ and thus
$$
\frac{|\alpha_3|}{\sigma^3\sqrt{n}} \leq \frac{\beta_3}{\sigma^3\sqrt{n}} 
= L_3 \leq L_4^{1/2} = \frac{1}{T_0}.
$$ 
Hence, according to definition (1.2), 
$|\Phi_3(x)| \leq c\,(1 + \frac{1}{T_0})$ for $x \leq 0$ and
$|1-\Phi_3(x)| \leq c\,(1 + \frac{1}{T_0})$ for $x \geq 0$, and thus
$|\Delta_n| \leq c\,(1 + \frac{1}{T_0})$.
This implies that (2.3) holds automatically in case $T_0 \leq 1$
for a suitable $c$. Thus, we may assume that $T_0 \geq 1$, i.e.,
$n \geq \beta_4/\sigma^4$.

In this case, the derivative of the function $G = \Phi_3$, which is given by
$$
\Phi_3'(x) = \varphi(x) + 
\frac{\alpha_3}{6\sigma^3\sqrt{n}}\,(x^3 - 3x)\, \varphi(x),
$$
is uniformly bounded in absolute value by some constant. Hence, 
by (2.1), for any $T_1 \geq T_0$,
\be
c \Delta_n \, \leq \, \int_0^{T_0} \frac{|f_n(t) - g_3(t)|}{t}\,dt +
\int_{T_0}^{T_1} \frac{|f_n(t) - g_3(t)|}{t}\,dt + \frac{1}{T_1}.
\en

It is known that $f_n(t)$ is approximated by $g_3(t)$ on the interval 
$|t| \leq 1/L_3$ with an error of order $1/n$ (using Taylor's expansion for 
$f(t)$ near zero and the product structure of $f_n(t)$). In particular, for 
a smaller interval $|t| \leq T_0$, there is a well-known estimate
$$
|f_n(t) - g_3(t)| \, \leq \,
c\,\frac{\beta_4}{\sigma^4 n}\,\min\{1,t^4\}\ e^{-t^2/8}
$$
(cf. e.g. [Bo2] for details). It allows one to properly bound the first integrand
in (2.4), which simplifies this Berry-Esseen estimate to the form
\be
c \Delta_n \leq \frac{\beta_4}{\sigma^4 n} + \frac{1}{T_1} + 
\int_{T_0}^{T_1} \frac{|f_n(t) - g_3(t)|}{t}\,dt.
\en

Now, according to (2.2) and using the assumption $T_0 \geq 1$, we also have 
\be
|g_3(t)| \leq \Big(1 + \frac{1}{6}\,t^3\Big)\,e^{-t^2/2} < 1.3\,e^{-t^2/8}
\qquad (t \geq 0),
\en
which implies
$$
\int_{T_0}^{T_1} \frac{|g_3(t)|}{t}\,dt \leq 
c\int_{T_0}^\infty e^{-t^2/8}\,dt < 4c\, e^{-T_0^2/8} < \frac{32\,c}{T_0^2} = 
32 c\, \frac{\beta_4}{\sigma^4 n}.
$$
As a result, (2.5) is simplified to
$$
c\,\Delta_n \leq \frac{\beta_4}{\sigma^4 n} + 
\frac{1}{T_1} + \int_{T_0}^{T_1} \frac{|f_n(t)|}{t}\,dt.
$$
Putting $T_1 = T\sigma \sqrt{n}$ and changing the variable, we arrive at 
(2.3). Note that the condition $T_1 \geq T_0$ 
is equivalent to $T \geq \frac{\sigma}{\sqrt{\beta_4}}$
\qed

\vskip5mm
Using Lemma 2.1, one obtains the statement of Theorem 1.1 in one direction.

\vskip5mm
{\bf Proposition 2.2.} {\sl Suppose that $\beta_4$ is finite and let, 
for some $p > 0$ and $q \in \R$,
$$
\frac{1}{1 - |f(t)|} =
O\Big(t^p\,(\log t)^q\Big) \quad {\sl as} \ \  t \rightarrow \infty.
$$
Then
\be
\Delta_n = 
O\Big(n^{-\frac{1}{2} - \frac{1}{p}}\,(\log n)^{\frac{q+1}{p}} + n^{-1}\Big).
\en
}

\vskip2mm
For $p < 2$ with arbitrary $q$ and for $p=2$ with $q \leq -1$, the relation (2.7) 
reduces to $\Delta_n = O\big(\frac{1}{n}\big)$, while in the other cases,
$$
\Delta_n = 
O\big(n^{-\frac{1}{2} - \frac{1}{p}}\,(\log n)^{\frac{q+1}{p}}\big).
$$
In particular, the hypothesis $\frac{1}{1 - |f(t)|} = \widetilde O(t^p)$ 
with $p \geq 2$ implies
$\Delta_n = \widetilde O(n^{-\frac{1}{2} - \frac{1}{p}})$.

\vskip5mm
{\bf Proof.} Suppose that $q \neq 0$. By the assumption, and since necessarily 
$X$ has a non-lattice distribution, we have
for all $T \geq t_0 = \frac{\sigma}{\sqrt{\beta_4}}$,
$$
M(T) = \max_{t_0 \leq t \leq T} |f(t)| \leq 
1 - \frac{a}{T^p\, \log^q(2+T)}
$$
with some constant $a > 0$. Using $1-u \leq e^{-u}$, we then get
$$
|f(t)|^n \leq M(T)^n \leq \exp\Big\{-\frac{na}{T^p\log^q(2 + T)}\Big\},
$$
so that
$$
\int_{t_0}^T \frac{|f(t)|^n}{t}\,dt \, \leq \,
\exp\Big\{-\frac{na}{T^p\log^q(2 + T)}\Big\}\,\log(T/t_0).
$$
Thus, by (2.3),
\be
c\,\Delta_n \leq \frac{\beta_4}{\sigma^4 n} + \frac{1}{T\sigma\sqrt{n}} + 
\exp\Big\{-\frac{na}{T^p\log^q(2 + T)}\Big\}\,\log(T/t_0).
\en

Let us take $T = T_n = (bn)^{1/p}\,(\log n)^{-r}$ with parameters
$r \geq 0$, $b>0$ to be precised later on and assuming that $n$ is
large enough. Then
$$
T_n^p \leq bn\,(\log n)^{-rp}, \qquad
\log(2 + T_n) \leq \frac{1}{p}\,\log n + O(\log \log n),
$$
and
$$
\log^q(2 + T_n) \leq 
\frac{1}{p^q}\,(\log n)^q + O\Big((\log n)^{q-1}\,\log \log n\Big).
$$
This gives 
$$
T_n^p\,\log^q(2+T_n) \leq \frac{b}{p^q}\,n\, (\log n)^{q-rp} + 
O\Big(n\,(\log n)^{q-rp-1}\,\log \log n\Big).
$$
Choosing $r = (q+1)/p$, the above is simplified to
$$
T_n^p\,\log^p(2+T_n) \leq \frac{b}{p^q}\,n\, (\log n)^{-1}\, 
\Big(1 + O\Big((\log n)^{-1}\,\log \log n\Big)\Big),
$$
and then
$$
\frac{na}{T_n^p\log^p(2 + T_n)} \geq 
\frac{ap^q}{b}\,\log n + O(\log \log n) \geq 2\log n,
$$
where the last inequality holds true with $b = ap^q/3$ for all $n$ large 
enough. In this case, the last term in (2.8) is estimated from above by 
$O(1/n)$.

In case $q=0$ with choice $r=1/p$, we clearly arrive at the same conclusion.
Therefore, (2.8) yields 
$$
\Delta_n = O\Big(\frac{1}{n} + \frac{1}{T_n\sqrt{n}}\Big) =
O\Big(\frac{1}{n} + n^{-\frac{1}{p} - \frac{1}{2}}\,(\log n)^r\Big), \qquad
r = \frac{q+1}{p}.
$$
\qed

\vskip10mm
\section{Non-uniform bounds based on uniform bounds}
\setcounter{equation}{0}

Suppose that a given distribution function $F$ is well approximated by
some function of bounded variation $G$ such that $G(-\infty) = 0$, 
$G(\infty) = 1$, in the sense of the Kolmogorov distance
$$
\Delta = \sup_x |F(x) - G(x)|.
$$
Based on this quantity, one would also like to see that $|F(x) - G(x)|$
decays polynomially fast for growing $x$. To this aim one may use moment 
assumptions together with some possible properties of $G$ related to its 
behavior at infinity.

\vskip5mm
{\bf Lemma 3.1.} {\sl Suppose that $F$ and $G$ have finite and equal second 
moments:
\be
\int_{-\infty}^\infty x^2\,dF(x) = \int_{-\infty}^\infty x^2\,dG(x).
\en
Then, for any $a>0$,
\begin{eqnarray}
\sup_x\, \Big[x^2\,|F(x)-G(x)|\Big]
 & \leq &
4a^2 \Delta + \int_{|x| \geq a} x^2\,dG(x) \nonumber \\
 & & + \ 
\max\Big\{\sup_{x \geq a}\, \big[x^2\,|1-G(x)|\big],
\sup_{x \leq -a} \big[x^2\,|G(x)|\big]\Big\}.
\end{eqnarray}
}

\vskip2mm
{\bf Proof.} For $|x| \leq a$, we have $x^2\,|F(x)-G(x)| \leq a^2 \Delta$ 
which is dominated by the right-hand side of (3.2). So, when estimating
$x^2\,|F(x)-G(x)|$, one may assume that $|x| > a$ and that
$\pm a$ are the points of continuity of both $F$ and $G$. 
Integrating by parts, we have
\bee
\int_{-a}^a y^2\,dF(y)
 & = &
a^2 (F(a)-G(a)) - a^2 (F(-a)-G(-a)) \\
 & & - \ 
2\int_{-a}^a y\,(F(y)-G(y))\,dy + \int_{-a}^a y^2\,dG(y). 
\ene
Hence
$$
\int_{-a}^a y^2\,dF(y) \geq -4a^2 \Delta + \int_{-a}^a y^2\,dG(y)
$$
which implies, by the moment assumption (3.1),
\be
\int_{|y| \geq a} y^2\,dF(y) \leq 4a^2 \Delta + \int_{|y| \geq a} y^2\,dG(y).
\en

On the other hand, in case $x \geq a$,
\bee
\int_{|y| \geq a} y^2\,dF(y) 
 & \geq &
\int_x^\infty y^2\,dF(y) \\
 & \geq &
x^2 (1-F(x)) \ = \ x^2 (G(x)-F(x)) + x^2\, (1-G(x)),
\ene
so,
$$
x^2 (G(x)-F(x)) \leq \int_{|y| \geq a} y^2\,dF(y) + 
\sup_{x \geq a} \big[x^2\, |1-G(x)|\big].
$$
Since also
$$
x^2 (F(x)-G(x)) \, \leq \, x^2 (1-G(x)) \, \leq \, \sup_{x \geq a} \big[x^2\, |1-G(x)|\big],
$$
we get
$$
x^2\,|F(x)-G(x)| \leq \int_{|y| \geq a} y^2\,dF(y) + 
\sup_{x \geq a}\, \big[x^2\,|1-G(x)|\big].
$$
By a similar argument, if $x \leq -a$,
$$
x^2\,|F(x)-G(x)| \leq \int_{|y| \geq a} y^2\,dF(y) + 
\sup_{x \leq -a} \big[x^2\,|G(x)|\big].
$$
Therefore, in both cases,
$$
x^2\,|F(x)-G(x)| \, \leq \, \int_{|y| \geq a} y^2\,dF(y) + 
\max\Big\{\sup_{x \geq a}\, \big[x^2\,|1-G(x)|\big], 
\sup_{x \leq -a} \big[x^2\,|G(x)|\big]\Big\}.
$$
It remains to involve (3.3).
\qed

\vskip5mm
In particular, if $G$ as measure is supported on the interval $[-a,a]$, then, 
under the moment assumption $(3.1)$, we have
\be
\sup_x\, \Big[x^2\,|F(x)-G(x)|\Big] \, \leq \, 4a^2 \Delta.
\en

In the general (non-compact) case, in order to optimize the inequality (3.2) 
over the variable $a$, an extra information is needed about the behavior 
of $G$. For example, let us require that, for some
parameters $A,B>0$,
\be
|G(x)| \leq A e^{-x^2/B} \ {\rm for} \ x \leq 0, \qquad
|1-G(x)| \leq A e^{-x^2/B} \ {\rm for} \ x \geq 0.
\en
The function $t e^{-t}$ is decreasing for $t \geq 1$. Hence,
if $x \geq a \geq \sqrt{B}$, we have
$$
x^2\,|1-G(x)| \leq A x^2\,e^{-x^2/B} \leq A a^2\,e^{-a^2/B}.
$$
In addition,
\bee
\int_a^\infty x^2\,dG(x)
 & = &
a^2\, (1-G(a)) + 2\int_a^\infty x\,(1-G(x))\,dx \\
 & \leq &
A a^2\,e^{-a^2/B} + 2A\int_a^\infty x\,e^{-x^2/B}\,dx
 \, = \,
A\, (a^2 + B)\,e^{-a^2/B} \, \leq \, 2A a^2\,e^{-a^2/B}.
\ene
Similar bounds also hold for the region $x \leq -a$.
Hence, the inequality (3.2) yields, for all $x \in \R$,
$$
x^2\,|F(x)-G(x)| \, \leq \,
4a^2 \Delta + 5A a^2\,e^{-a^2/B}, \qquad a \geq \sqrt{B}.
$$
Moreover, choosing
$a^2 = B\,\log(e + \frac{1}{\Delta})$, the above right-hand side
becomes
$$
4B\,\Delta\log\Big(e + \frac{1}{\Delta}\Big) + 
5AB\,\frac{1}{e + \frac{1}{\Delta}}\log\Big(e + \frac{1}{\Delta}\Big) \, \leq \,
(4B + 5AB)\,\Delta\log\Big(e + \frac{1}{\Delta}\Big).
$$

Note that the parameters $A$ and $B$ may not be arbitrary.
Applying the hypothesis (3.5) at the origin $x=0$, we get
$
1 \leq |G(0)| + |1-G(0)| \leq 2A.
$
So, necessarily $A \geq \frac{1}{2}$ and hence $4+5A \leq 13 A$.
Thus, applying Lemma 3.1, we arrive at the following assertion.

\vskip5mm
{\bf Proposition 3.2.} {\sl Under the assumptions $(3.1)$ and $(3.5)$,
\be
\sup_x\,\Big[x^2\,|F(x)-G(x)|\Big] \, \leq \,
13\,AB\, \Delta \log\Big(e + \frac{1}{\Delta}\Big).
\en
}

\vskip2mm
In case of the normal distribution function $G = \Phi$, we have
$1 - \Phi(x) \leq \frac{1}{2}\,e^{-x^2/2}$ ($x \geq 0$),
so, the conditions (3.1) and (3.5) are fulfilled with
$A = \frac{1}{2}$ and $B = 2$. Hence
\be
\sup_x\,\Big[x^2\,|F(x)-\Phi(x)|\Big] \, \leq \,
13\, \Delta \log(e + 1/\Delta),
\en
provided that $\int_{-\infty}^\infty x^2\,dF(x) = 1$. In fact, this bound can be
generalized in order to control a polynomial decay of $|F(x)-\Phi(x)|$ of any order $p>0$. 
Namely, if $\Delta \leq \frac{1}{\sqrt{e}}$, one has
$$
\sup_x \Big[(1 + |x|^p)\,|F(x)-\Phi(x)|\Big] \, \leq \,
C_p\, \Delta \log^{p/2}(1/\Delta) + \lambda_p,
$$
where
$$
\lambda_p = 
\bigg|\int_{-\infty}^\infty |x|^p\,dF(x) - \int_{-\infty}^\infty |x|^p\,d\Phi(x)\bigg|,
$$
and the constant $C_p$ depends on $p$ only. This inequality can be found in 
[P2], Ch.\,V, Theorem\,11, pp. 174-176 
(where it is attributed to Kolodyazhnyi [K]). The proof of Lemma 3.1 given above 
follows the same line of arguments as in [P2].
As for Proposition 3.2, we will need with $G = \Phi_3$.

\vskip10mm
\section{Deviations of characteristic functions}
\setcounter{equation}{0}

The non-uniform bound (3.6) allows one to control deviations of
the Fourier-Stieltjes transform $f$ of the distribution function $F$ from 
the Fourier-Stieltjes transform of $G$. Recall that $G$ is assumed to be a function 
of bounded variation such that $G(-\infty) = 0$ and $G(\infty) = 1$.

From (3.6) it follows that, for any $b>0$,
$$
\sup_x\,\Big[(b^2 + x^2)\,|F(x)-G(x)|\Big] \, \leq \,
b^2 \Delta + 13\,AB\, \Delta \log\Big(e + \frac{1}{\Delta}\Big),
$$
and therefore
\bee
W_1(F,G)
 & \equiv &
\int_{-\infty}^\infty |F(x) - G(x)|\,dx \\
 & \leq &
\frac{\pi}{b}\,
\Big[b^2 \Delta + 13\,AB\, \Delta \log\Big(e + \frac{1}{\Delta}\Big)\Big] 
 \, = \,
\pi \Delta\,
\Big[b + \frac{13\,AB}{b}\, \log\Big(e + \frac{1}{\Delta}\Big)\Big].
\ene
Optimizing the right-hand side over all $b>0$ and using 
$\pi \sqrt{26} < 16.02$, we arrive at
\be
W_1(F,G) \, \leq \,
16.02\,\sqrt{AB}\ \Delta \log^{1/2}\Big(e + \frac{1}{\Delta}\Big).
\en
In particular, we get:

\vskip5mm
{\bf Proposition 4.1.} {\sl Under the assumptions $(3.1)$ and $(3.5)$,
for all $t \in \R$,
\be
|f(t) - g(t)| \, \leq \, 16.02\,\sqrt{AB}\ |t|\, 
\Delta \log^{1/2}\Big(e + \frac{1}{\Delta}\Big),
\en
where $\Delta = \sup_x |F(x) - G(x)|$. 
}

\vskip5mm
This bound follows from (4.1) via the the identity
$$
f(t) - g(t) = -it \int_{-\infty}^\infty e^{itx}\,(F(x)-G(x))\,dx.
$$

The logarithmic term in (4.2) may be removed for compactly supported 
distributions $G$, even if $F$ is not compactly supported.
Indeed, starting from (3.4), for any $b>0$,
$$
\sup_x\,\Big[(b^2 + x^2)\,|F(x)-G(x)|\Big] \, \leq \, 
b^2 \Delta + 4a^2\, \Delta,
$$
and therefore
$$
W_1(F,G) \, \leq \, \frac{\pi}{b}\,(b^2 + 4a^2)\, \Delta \, = \, 
\pi \Delta\, \Big[b + \frac{4a^2}{b}\Big] = 4\pi a\,\Delta,
$$
where in the last equality we take an optimal value $b = 2a$. Hence,
if $G$ is supported on the interval $[-a,a]$ (as measure) and has 
the same second moment as $F$, then
$$
|f(t) - g(t)| \, \leq \, 4\pi a\, \Delta\,|t| \qquad (t \in \R).
$$

\vskip2mm
Now, let us return to the setting of Theorem 1.1 and specialize
Proposition 4.1 to
$$
G(x) = \Phi_3(x) = \Phi(x) - 
\frac{\alpha_3}{6\sigma^3\sqrt{n}}\,(x^2-1)\, \varphi(x).
$$ 
As was explained in Section 2,
$\frac{|\alpha_3|}{\sigma^3\sqrt{n}} \leq 1$ as long as
$n \geq \beta_4/\sigma^4$. In this case, for any $x \geq 0$,
$$
|1 - \Phi_3(x)| \, \leq \,
|1 - \Phi(x)| + \frac{|\alpha_3|}{6\sigma^3\sqrt{n}}\,|x^2-1|\, \varphi(x)
 \, \leq \, \frac{1}{2}\,e^{-x^2/2} + 
\frac{1}{6\sqrt{2\pi}}\, |x^2-1|\, e^{-x^2/2}.
$$
Being multiplied by $e^{x^2/4}$, the above right-hand side attains maximum 
at zero, hence
$$
|1 - \Phi_3(x)| \, \leq \,
\Big(\frac{1}{2} + \frac{1}{6\sqrt{2\pi}}\Big)\,e^{-x^2/4} \, < \, 
0.57\,e^{-x^2/4}.
$$
Thus, the assumption (3.5) is fulfilled with $A = 0.57$ and $B = 4$.
We then get:

\vskip5mm
{\bf Corollary 4.2.} {\sl Suppose that $\beta_4$ is finite. For all 
$n \geq \beta_4/\sigma^4$, the characteristic function $f_n(t)$ of $Z_n$
satisfies, for all $t \in \R$,
\be
|f_n(t) - g_3(t)| \, \leq \, 24.2\, |t|\, 
\Delta_n \log^{1/2}\Big(e + \frac{1}{\Delta_n}\Big),
\en
where $\Delta_n = \sup_x |F_n(x) - \Phi_3(x)|$.
}

\vskip5mm
In fact, when $\alpha_3 = 0$, we have $\Phi_3 = \Phi$,
and the requirement $n \geq \beta_4/\sigma^4$ together with the 4-th moment
assumption are not needed in Corollary 4.2. Moreover, since $AB = 1$ for 
$G = \Phi$ in (3.5), from (3.7) we obtain a better numerical constant. Namely,
$$
|f_n(t) - e^{-t^2/2}| \, \leq \, 16.02\, |t|\, 
\Delta_n \log^{1/2}\Big(e + \frac{1}{\Delta_n}\Big).
$$

\vskip10mm
\section{Necessity part in Theorem 1.1}
\setcounter{equation}{0}

Keeping the setting of Theorem 1.1,
one may use the deviation inequality (4.3) to show that $f(t)$
is properly bounded away from 1 and thus to reverse the statement
of Proposition 2.2. In this direction, only the finiteness of the 3-rd
absolute moments is needed (which is necessary, since $\alpha_3$ 
participates in the definition of $\Phi_3$).

\vskip5mm
{\bf Proposition 5.1.} {\sl Suppose that, for some $p>0$ and
$q \in \R$,
$$
\Delta_n = O\Big(n^{-(\frac{1}{2} + \frac{1}{p})}\,(\log n)^q\Big) \ 
\quad {\sl as} \ \  n \rightarrow \infty.
$$
Then
\be
\frac{1}{1 - |f(t)|} =
O\Big(t^p\,(\log t)^{p\,(\frac{1}{2} + q)}\Big) \quad {\sl as} \ \  
t \rightarrow \infty.
\en
}

\vskip5mm
{\bf Proof.} By the assumption,
$$
\Delta_n \log^{1/2}\Big(e + \frac{1}{\Delta_n}\Big) =
O\Big(n^{-\frac{1}{2} - \frac{1}{p}}\,(\log n)^{q + \frac{1}{2}}\Big).
$$
Hence, using the upper bound (2.6) on $|g_3(t)|$, (4.3) yields, for all 
$n \geq \beta_4/\sigma^4$,
$$
|f_n(t)| \, \leq \, 1.3\, e^{-t^2/8} + 
c\, |t|\,n^{-\frac{1}{2} - \frac{1}{p}}\,(\log n)^{q + \frac{1}{2}}.
$$
Here in the region $t \geq \sqrt{n}$, the second term on the right-hand side
dominates the first one. Replacing $t$ with $t\sqrt{n}$, we therefore
obtain that
\be
|f(t/\sigma)|^n \leq c_{p,q}\, t\,n^{-1/p}\,(\log(n+1))^{q + 1/2}, \qquad 
t \geq 1,
\en
with some $(p,q)$-dependent constant $c_{p,q}$. Assuming that $t \geq e$,
let us choose 
\be
n = [2A t^p\, (\log t)^r]
\en
with parameters $A \geq 1$ and $r>0$. In this case,
$$
n^{-1/p} \leq \big(A t^p\, (\log t)^r\big)^{-1/p} =
A^{-1/p}\,t^{-1}\, (\log t)^{-r/p}
$$
and
\bee
\log(n+1) 
 & \leq &
\log(4A t^p\, (\log t)^r) \, = \, \log(4A) + p \log t + r \log \log t \\
 & < &
\log(4A) + (p+r) \log t \, < \, (p+r+1) \log t,
\ene
where in the last inequality we require that $t \geq 4A$.
Hence
\bee
n^{-1/p}\,(\log(n+1))^{q + 1/2} 
 & \leq &
A^{-1/p}\,t^{-1}\, (\log t)^{-r/p} \cdot
(p+r+1)^{q + 1/2}\, (\log t)^{q + 1/2} \\ 
 & = &
A^{-1/p}\,c'_{p,q} \, t^{-1},
\ene
where we chose $r = p (q + 1/2)$ on the last step.
Hence, with some $(p,q)$-dependent constant, (5.2) is simplified to
$$
|f(t/\sigma)|^n \leq c_{p,q}\,A^{-1/p},
$$ 
which can be made smaller than $1/e$ by choosing a sufficiently large
value of $A$. Thus, recalling (5.3), we have
$$
|f(t/\sigma)| \leq e^{-1/n} \leq 1 - \frac{1}{2n} \leq 
1 - \frac{1}{4A t^p\, (\log t)^r},
$$
which yields (5.1).
\qed

\vskip10mm
\section{Diophantine inequalities}
\setcounter{equation}{0}

Turning to Corollary 1.2 and other applications of Theorem 1.1, it makes 
sense to describe a somewhat more general situation. First let us list a few 
simple metric properties of the function $x \rightarrow \|x\|$ in the real variable $x$.
This function is even, 1-periodic, and satisfies, for all real $x,y$,

\vskip2mm
(i) \ $\|x\| \leq |x|$;

(ii) \ $\|x + y\| \leq \|x\| + \|y\|$;

(iii)\ $|\,\|x\| - \|y\|\,| \leq \|x - y\|$.

\vskip2mm
\noindent
In addition,
\be
|\cos(\pi x)| \leq \exp\{-\pi^2 \|x\|^2/2\}, \qquad
4\, \|x\|^2 \leq 1 - |\cos(\pi x)| \leq \frac{\pi^2}{2}\,\|x\|^2.
\en
The inequalities in (6.1) are elementary, and we omit the proofs.

Below, we denote by $n(x)$ the closest integer to $x$, so that
$\|x\| = |x - n(x)|$ (for definiteness, let $n(x) = n$ in case $x = n + 1/2$).

\vskip5mm
{\bf Lemma 6.1.} {\sl Given real numbers $\alpha_1,\dots,\alpha_m$, 
suppose that $\max_{k \leq m} \|n\alpha_k\| \geq \ep(n) > 0$ 
for all integers $n \geq 1$. Then, for all $t \geq 1$~real,
\be
\|t\|^2 + \|t\alpha_1\|^2 + \dots + \|t\alpha_m\|^2 \geq c^2 \ep(n(t))^2,
\en
where $c^{-1} = 1 + \max_{k\leq m} |\alpha_k|$.
}

\vskip5mm
{\bf Proof.} One may assume that all $\alpha_k > 0$.
Let $t = n + \gamma$, $|\gamma| = \|t\|$, with $n = n(t)$.
If $\|t\| \geq c\ep(n)$, $c>0$, then automatically 
$$
M(t) \equiv 
\max\big\{\|t\|,\|t\alpha_1\|,\dots,\|t\alpha_m\|\big\} \geq c \ep(n).
$$
Now, suppose that $\|t\| < c\ep(n)$. By the assumption, $\|n\alpha_k\| \geq \ep(n)$ 
for some $k \leq m$. Since $t\alpha_k = n\alpha_k + \gamma \alpha_k$, 
we get, applying the properties (i) and (iii):
\bee
\|t\alpha_k\| 
 & \geq &
\|n \alpha_k\| - \|\gamma \alpha_k\| \\
 & \geq &
\|n \alpha_k\| - |\gamma \alpha_k| \, = \,
\|n \alpha_k\| - \|t\|\, \alpha_k \, \geq \, (1 - c\alpha_k)\, \ep(n).
\ene
Here $1 - c\alpha_k = c$ for $c = \frac{1}{1 + \alpha_k}$, and then
$\|t\alpha_k\| \geq c\ep(n)$ in both cases. Hence, 
$M(t) \geq \frac{\ep(n)}{1 + |\alpha_k|}$.
\qed

\vskip5mm
Clearly, (6.2) with integer values $t=n$ returns us to
the assumption, up to an $\alpha_k$-depending factor in front of $\ep(n)$.

Let us now consider a system of $m$ Diophantine inequalities
$$
\Big|\alpha_k - \frac{r_k}{n}\Big| < \frac{\ep(n)}{n}, \qquad
k = 1,\dots,m \quad (n \geq 1),
$$
about which one is usually concerned whether or not it has infinitely
many integer solutions $(r_1,\dots,r_m,n)$. Here,
we choose the particular functions 
$\ep(n) = c\,n^{-\eta}\,(\log(n+1))^{-\eta'}$ and consider the opposite property:
\be
\liminf_{n \rightarrow \infty} \Big[\,n^\eta\,(\log n)^{\eta'}\,
\max\{\|n\alpha_1\|,\dots,\|n\alpha_m\|\}\Big] > 0.
\en
One may rephrase this in terms of the characteristic function
\be
f(t) = \cos(t)\,\cos(\alpha_1 t) \dots \cos(\alpha_m t)
\en
of the sum $X = \xi_0 + \alpha_1 \xi_1 + \dots + \alpha_m \xi_m$, 
where $\xi_k$ are independent Bernoulli random variables, taking the values 
$\pm 1$ with probability $1/2$.

\vskip5mm
{\bf Lemma 6.2.} {\sl Given $\alpha_1,\dots,\alpha_m \in \R$ and $\eta>0$, 
$\eta' \in \R$, the relation $(6.3)$ is equivalent to the property that 
the characteristic function $f$ in $(6.4)$ satisfies
\be
\frac{1}{1-|f(t)|} = O\big(t^{2\eta}\,(\log t)^{2\eta'}\big) \quad
{\sl as} \ \ t \rightarrow \infty.
\en
}

\vskip2mm
{\bf Proof.} For (6.3) to hold, it is necessary that at least one of 
$\alpha_k$ be irrational. Moreover, this relation may be strengthened to
\be
\max_{1 \leq k \leq m}\, \|n\alpha_k\| \, \geq \, 
\frac{c}{n^\eta\,(\log(n+1))^{\eta'}}, \qquad n \geq 1,
\en
with some constant $c>0$ independent of $n$. Moreover, according to Lemma 6.1 with
$\ep(n)$ as above, we see that (6.6) is equivalent to
\be
\|t\|^2 + \|t\alpha_1\|^2 + \dots + \|t\alpha_m\|^2 \, \geq \,
\frac{c}{t^{2\eta}\,(\log t)^{2\eta'}}, \quad t \geq 2 \ ({\rm real})
\en
(modulo positive constants). Combining (6.7) with the first inequality 
in (6.1) yields
$$
|f(\pi t)| \, \leq \,
\exp\Big\{-\frac{\pi^2}{2}\, 
\big(\|t\|^2 + \|\alpha_1 t\|^2 + \dots + \|t\alpha_m\|^2\big)\Big\}\, \leq \, 
\exp\Big\{-\frac{c}{t^{2\eta}\,(\log t)^{2\eta'}}\Big\},
$$
which thus leads to the required relation (6.5).

Conversely, (6.5) yields
\be
1 - |f(\pi t)| \, \geq \, \frac{c}{t^\eta\,(\log(t+1))^{\eta'}}, \qquad t \geq 1,
\en
so that for the integer values $t = n$ we get
$$
1 - \frac{c}{n^{2\eta}\,(\log(n+1))^{2\eta'}} \, \geq \, |f(\pi n)| =
(1-\delta_1)\dots (1 - \delta_m), \qquad \delta_k = 1-|\cos(\pi n \alpha_k)|.
$$
Since the right-hand side is greater than or equal to
$1 - (\delta_1 + \dots + \delta_m)$, we obtain
$$
\frac{c}{n^{2\eta}\,(\log(n+1))^{2\eta'}} \leq \delta_1 + \dots + \delta_m.
$$
Recalling (6.1), we have 
$\delta_k \leq \frac{\pi^2}{2}\,\|n \alpha_k\|^2$ and thus
$$
\frac{c}{n^{2\eta}\,(\log(n+1))^{2\eta'}} \, \leq \,
\frac{\pi^2}{2}\, \sum_{k=1}^m \|n \alpha_k\|^2 \, \leq \,
\frac{m \pi^2}{2}\, \max_{k \leq m} \|n \alpha_k\|^2.
$$
This gives (6.6) and therefore (6.3).
\qed

\vskip5mm
A similar conclusion continues to hold for other
characteristic functions including
\be
f(t) = p_0 \cos(t) + \sum_{k=1}^m\,p_k \cos(\alpha_k t),
\en
where $p_k$ are fixed positive parameters such that $p_0 + \dots + p_m = 1$.
Indeed, by (6.1),
\bee
|f(\pi t)| 
 & \leq & 
p_0\,\big(1 - 4\, \|t\|^2) + 
\sum_{k=1}^m p_k \big(1 - 4\, \|\alpha_k t\|^2\big) \\
 & \leq &
1 - p' \Big(\|t\|^2 + \|\alpha_1 t\|^2 + \dots + \|\alpha_m t\|^2\Big), \qquad
p' = 4\,\min_{0 \leq k \leq m} p_k.
\ene
Starting from (6.6)-(6.7), we would obtain again (6.5).

Conversely, (6.5) leads to (6.8), which at the even integer values $t = 2n$ 
yields
$$
1 - \frac{c}{n^{2\eta}\,(\log(n+1))^{2\eta'}} \, \geq \,
f(2\pi n) \, = \, p_0 + \sum_{k=1}^m p_k \cos(2\pi n \alpha_k) \, = \, 
1 - 2\sum_{k=1}^m p_k \delta_k^2,
$$
where now $\delta_k = \sin(\pi n\alpha_k)$. Using 
$|\sin(\pi x)| \leq \pi\,\|x\|$, the above inequality yields
$$
\frac{c}{n^{2\eta}\,(\log(n+1))^{2\eta'}} \, \leq \,
2\pi^2 \sum_{k=1}^m p_k \|n \alpha_k\|^2\, \leq \, 
2\pi^2\, \max_{k \leq m} \|n \alpha_k\|^2.
$$
As a result, we arrive at:

\vskip5mm
{\bf Lemma 6.3.} {\sl The assertion of Lemma $6.2$ is also true
for all characteristic functions $f$ of the form $(6.9)$.
}

\vskip5mm
We are prepared to prove Corollary 1.2, in fact -- in a more
precise and general form, if we apply Propositions 2.2 and 5.1. 
Let us return to the setting of Theorem 1.1 in which we will assume
that the random variable $X$ has a characteristic function $f$ given by 
(6.4) or (6.9). Equivalently, if we denote by
$B_\alpha = \frac{1}{2}\,\delta_{\alpha} + \frac{1}{2}\,\delta_{-\alpha}$
the symmetric Bernoulli measure supported on $\{-\alpha,\alpha\}$,
the distribution $F$ of $X$ may be written
(as measure) in either of the two forms
$$
F = B_1 * B_{\alpha_1} * \dots * B_{\alpha_m}, \qquad
F = p_0 B_1 + \sum_{k=1}^m p_k B_{\alpha_k} \quad (p_k > 0, \
p_0 + \dots + p_m = 1).
$$
Since any such measure is symmetric about the origin, the uniform distance 
in Theorem 1.1 is defined by
$
\Delta_n = \sup_x |F_n(x) - \Phi(x)|.
$

\vskip5mm
{\bf Proposition 6.4.} {\sl Given $\alpha_1,\dots,\alpha_m \in \R$, 
suppose that with some $\eta \geq 1$, $\eta' \in \R$, 
\be
\liminf_{n \rightarrow \infty} \Big[\,n^\eta\,(\log n)^{\eta'}\,
\max\big\{\|n\alpha_1\|,\dots,\|n\alpha_m\|\big\}\Big] > 0.
\en
Then
\be
\Delta_n = 
O\Big(n^{-\frac{1}{2} - \frac{1}{2\eta}}\,(\log n)^{\eta''}\Big)
\en
with $\eta'' = \frac{2\eta'+1}{2\eta}$ in case $\eta > 1$ and
$\eta'' = \max\big\{\frac{2\eta'+1}{2},0\big\}$ in case $\eta = 1$. 

\vskip2mm
\noindent
Conversely, if $(6.11)$
holds with some $\eta>0$, $\eta'' \in \R$, then $(6.10)$ is
fulfilled with $\eta' = \eta\,(\frac{1}{2} + \eta'')$.
}

\vskip5mm
Indeed, starting from the hypothesis (6.10), we obtain (6.5), so that 
the condition of Proposition 2.2 is fulfilled with $p = 2\eta$ 
and $q = 2\eta'$. Hence, by Proposition 2.2,
$$
\Delta_n = 
O\Big(n^{-\frac{1}{2} - \frac{1}{p}}\,(\log n)^{\frac{q+1}{p}} + n^{-1}\Big),
$$
i.e. (6.11). Conversely, (6.11) ensures that the condition of Proposition 5.1 
is fulfilled with $p = 2\eta$ and $q = \eta''$. Therefore,
$$
\frac{1}{1 - |f(t)|} =
O\Big(t^p\,(\log t)^{p\,(\frac{1}{2} + q)}\Big) =
O\Big(t^{2\eta}\,(\log t)^{2\eta\,(\frac{1}{2} + \eta'')}\Big),
$$
which is (6.5) with $2\eta' = 2\eta\,(\frac{1}{2} + \eta'')$.

\vskip10mm
\section{Special values of $\alpha$ and typical behavior of $\Delta_n$}
\setcounter{equation}{0}

Let us restrict the setting of Proposition 6.4 to the case $m=1$ 
and assume that the distribution $F$ of $X$ has a convolution structure, i.e.,
$X = X' + \alpha X''$, where $X',X''$ are independent random variables with 
a symmetric Bernoulli distribution on $\{-1,1\}$. The corresponding 
characteristic function is then given by $f(t) = \cos(t)\,\cos(\alpha t)$,
and the second moment of $F$ is $\sigma^2 = 1 + \alpha^2$. Hence, 
the measure $F_n$ from Theorem 1.1 represents the distribution of
$$
Z_n = 
\frac{1}{\sqrt{1 + \alpha^2}}\ Z_n' + \frac{\alpha}{\sqrt{1 + \alpha^2}}\ Z_n'',
$$
where $Z'$ and $Z_n''$ are independent normalized sums of $n$
independent copies of $X'$ and $X''$. 

Put 
$$
\Delta_n(\alpha) = \sup_x\, |F_n(x) - \Phi(x)|.
$$
Since $\P\{Z_n = 0\} \geq \P\{Z_n' = 0\} \, \P\{Z_n'' = 0\} > \frac{c}{n}$,
we necessarily have
$\Delta_n(\alpha) > \frac{c}{n}$ with some absolute constant $c>0$.
On the other hand, Proposition 6.4 implies:

\vskip5mm
{\bf Corollary 7.1.} {\sl If
\be
\liminf_{n \rightarrow \infty} \Big[\,n^\eta\,(\log n)^{\eta'}\,
\|n\alpha\|\Big] > 0,
\en
for some $\eta \geq 1$, $\eta' \in \R$, then
\be
\Delta_n(\alpha) = 
O\Big(n^{-\frac{1}{2} - \frac{1}{2\eta}}\,(\log n)^{\eta''}\Big)
\en
with $\eta'' = \frac{2\eta'+1}{2\eta}$. In turn, the latter relation
implies $(7.1)$ with $\eta' = \eta\,(\frac{1}{2} + \eta'')$.
}

\vskip5mm
This is a more precise formulation of Corollary 1.2.
Note that (7.1) is impossible for $\eta = 1$ and $\eta' < 0$
(by Dirichlet's theorem), so that necessarily
$\eta'' =  \max\big\{\frac{2\eta'+1}{2},0\big\} = \frac{2\eta'+1}{2} \geq 
\frac{1}{2}$. Similarly, (7.2) is impossible for $\eta = 1$ and $\eta'' < 0$.

The relation (7.1) with $\eta=1$, $\eta' = 0$ defines the class 
of the so-called badly approximable numbers $\alpha$ which can be characterized
in terms of continued fractions. Namely, representing
$$
\alpha = a_0 + \frac{1}{a_1 + \frac{1}{a_2 + \frac{1}{a_3 + \ldots}}}\, ,
$$
where $a_0$ is an integer and $a_1,a_2,\dots$ are positive integers,
the property of being badly approximable is equivalent to
$\sup_i a_i < \infty$. In particular, all quadratic irrationalities
(e.g. $\alpha = \sqrt{2}$) belong to this class, cf. [S2]. 
Since in this case $\eta''  = \frac{2\eta'+1}{2} = \frac{1}{2}$, we arrive at:

\vskip5mm
{\bf Corollary 7.2.} {\sl For any badly approximable number $\alpha$, we have
$
\Delta_n(\alpha) = O\big(\frac{1}{n}\sqrt{\log n}\,\big).
$
}

\vskip5mm
It is not clear at all whether one can improve this rate for at least
one $\alpha$. On the other hand, at the expense of a logarithmic term,
one may involve almost all values of $\alpha$. To this aim, one may
apply a theorem due to Khinchine which asserts the following 
(cf. [C], [S2]).
Suppose that a function $\psi(n)>0$ is defined on the positive
integers. If $\psi(n)$ is non-increasing and $\sum_{n=1}^\infty \psi(n) = \infty$, then the inequality
\be
\Big|\alpha - \frac{p}{n}\Big| < \frac{\psi(n)}{n}
\en
has infinitely many integer solutions $(p,n)$ for almost all $\alpha$ 
(with respect to the Lebesgue measure on the real line). But when
$\sum_{n=1}^\infty \psi(n) < \infty$, 
(7.3) has only finitely many solutions for almost all $\alpha$.
This second assertion is an easy part of Khinchine's theorem, which
may be quantified in terms of the function
$$
r_\psi(\alpha) = \inf_{n \geq 1} \Big[\,\frac{1}{\psi(n)}\,\|n\alpha\|\Big].
$$
Indeed, restricting ourselves 
(without loss of generality) to the values $0 < \alpha < 1$,
first note that, for any integer $n \geq 1$ and $\delta > 0$,
$$
{\rm mes}\{\alpha \in (0,1): \|n\alpha\| < \delta\} \, \leq \, 2\delta
$$
(with equality in case $\delta \leq 1/2$). Hence, for any $r>0$,
$$
{\rm mes}\{\alpha \in (0,1):r_\psi(\alpha) < r\} \, \leq \, \sum_{n=1}^\infty 
{\rm mes}\Big\{\alpha \in (0,1):\frac{1}{\psi(n)}\,\|n\alpha\| < r\Big\} \, \leq \, 
\sum_{n=1}^\infty 2r\, \psi(n),
$$
and thus
$$
{\rm mes}\{\alpha \in (0,1):r_\psi(\alpha) < r\} \, \leq \, Cr \qquad
(r>0)
$$
with constant $C = 2\sum_{n=1}^\infty \psi(n)$. In particular,
$r_\psi(\alpha) > 0$ for almost all $\alpha$.

For example, choosing the sequence $\psi(n) = 1/(n\log^{1+\ep}(n+1))$, 
Corollary 7.1 provides a rate which is applicable to almost all $\alpha$.

\vskip5mm
{\bf Corollary 7.3.} {\sl Given $\ep>0$, for almost all $\alpha \in \R$, 
we have
$
\Delta_n(\alpha) = O\big(\frac{1}{n}\,(\log n)^{3/2 + \ep}\big).
$
}

\vskip5mm
It is not clear whether or not the power of the logarithmic term may be improved.
At least, this is possible on average when $\alpha$ varies inside a given interval, 
say $0 < \alpha < 1$.

\vskip5mm
{\bf Proposition 7.4.} {\sl With some absolute constant $c>0$,
for all $n \geq 1$,
\be
\int_0^1 \Delta_n(\alpha)\,d\alpha \, \leq \, c\,\frac{\log(n+1)}{n}.
\en
}

{\bf Proof.}
Our basic tool is the Berry-Esseen inequality of Lemma 2.1.
For the distribution $F$, we have $\alpha_3 = \E X^3 = 0$ and
$$
\beta_4 = \E X^4 = \E\, (X' + \alpha X'')^4 = 1 + 6\alpha^2 + \alpha^4.
$$
In order to control the integral in (2.3), recall that $\sigma^2 = 1 + \alpha^2$ and 
note that $\sigma^4 \leq \beta_4 \leq 2\sigma^4$. Using 
$\frac{\sigma}{\sqrt{\beta_4}} \geq \frac{1}{\sqrt{2(1 + \alpha^2)}} \geq 
\frac{1}{2}$, Lemma 2.1 with $T = \sqrt{n}$ gives that
\be
c\,\Delta_n(\alpha) \leq \frac{1}{n} + I_n(\alpha), \quad {\rm where} \quad 
I_n(\alpha) = \int_{1/2}^{\sqrt{n}} \frac{|\cos(t)\,\cos(\alpha t)|^n}{t}\,dt.
\en
By simple calculus, for any $t \geq 1/2$,
$$
\psi_n(t) \equiv \int_0^t |\cos(s)|^n\,ds \leq \frac{t}{\sqrt{n}}\sqrt{2\pi},
$$
so
$$
\int_0^1 I_n(\alpha)\,d\alpha  \, = \,
\int_{1/2}^{\sqrt{n}} \frac{|\cos t|^n}{t^2}\ \psi_n(t)\,dt \, \leq \,
\frac{\sqrt{2\pi}}{\sqrt{n}}\, \int_{1/2}^{\sqrt{n}} \frac{|\cos t|^n}{t}\,dt \, \leq \,
\frac{c\,\log(n+1)}{n}.
$$
Thus, integrating the inequality in (7.5) over $\alpha$, we are led to (7.4).
\qed

\vskip5mm
{\bf Remarks.}
Corollary 7.3 with quantity $\Delta_n(\alpha) = \sup_x\, |F_n(x) - \Phi_3(x)|$
remains to hold in a more general situation
$X = X' + \alpha X''$, where $X',X''$ are independent random variables with 
non-degenerate distributions and finite 4-th absolute moments.
This extension requires an extra analysis of the behavior
of characteristic functions, and we will discuss it somewhere else.
Let us note that it is possible to improve the rate of convergence 
(in particular, to remove the logarithmic term) in models such as
$X = X^{(0)} + \alpha_1 X^{(1)} + \dots + \alpha_m X^{(m)}$ with
$m \geq 2$ independent summands $X^{(k)}$. See also [K-S] on randomized
versions of the central limit theorem.

\end{document}